\documentclass[notitlepage]{amsart}

\usepackage{amsfonts}
\usepackage[latin1]{inputenc}
\usepackage[all]{xypic}
\newtheorem{pro}{Proposition}[section]
\newtheorem{teo}[pro]{Theorem}
\newtheorem{defi}[pro]{Definition}
\newtheorem{lem}[pro]{Lemma}
\newtheorem{cor}[pro]{Corollary}
\newtheorem{rk}[pro]{Remark}
\newtheorem{ex}[pro]{Example}
\newtheorem{exs}[pro]{Examples}

\newcommand{\pd}{{\mathrm{pd}}}
\newcommand{\mini}{{\mathrm{min}}}
\newcommand{\maxi}{{\mathrm{max}}}
\newcommand{\modu}{{\mathrm{mod}}}
\newcommand{\Ima}{{\mathrm{Im}}}
\newcommand{\Ker}{{\mathrm{Ker}}}

\newcommand{\Trace}{{\mathrm{Trace}}}
\newcommand{\Id}{{1_\mathcal{C}}}

\newcommand{\Hom}{{\mathrm{Hom}}}
\newcommand{\F}{{\mathcal{F}}}
\newcommand{\Z}{{\mathbb{Z}}}

\newcommand{\T}{{\mathcal{T}}}
\newcommand{\C}{{\mathcal{C}}}
\newcommand{\X}{{\mathcal{X}}}
\newcommand{\Y}{{\mathcal{Y}}}
\newcommand{\add}{{\mathrm{add}}}
\newcommand{\rad}{{\mathrm{rad}}}
\newcommand{\soc}{{\mathrm{soc}}}
\newcommand{\tops}{{\mathrm{top}}}

\newcommand{\fun}{{\mathrm{End}_{\Z}\,(\C)}}
\newcommand{\ab}{{{(\alpha,\beta)}}}
\newcommand{\dl}{{{\ell\ell}}}

\newcommand{\Fa}{{\mathrm{F}_{\alpha}}}
\newcommand{\Ga}{{\mathrm{G}_{\alpha}}}

\newcommand{\FFa}{{\mathcal{F}_{\alpha}}}
\newcommand{\Ta}{{\mathcal{T}_{\alpha}}}
\newcommand{\qa}{{q_{\alpha}}}
\newcommand{\qb}{{q_{\beta}}}
\newcommand{\qtp}{{q_{t'}}}
\newcommand{\FFb}{{\mathcal{F}_{\beta}}}
\newcommand{\Tb}{{\mathcal{T}_{\beta}}}
\newcommand{\Tq}{{\mathcal{T}_{\qb}}}
\newcommand{\FFq}{{\mathcal{F}_{\qb}}}
\newcommand{\FFqa}{{\mathcal{F}_{\qa}}}

\newcommand{\Gq}{{\emph{G}_{\qa}}}

\newcommand{\tp}{{(\mathcal{T},\mathcal{F})}}
\newcommand{\tpp}{{(\mathcal{T'},\mathcal{F'})}}
\newenvironment{dem}{\noindent {\bf Proof.}}{\hfill $\Box$\\}

\usepackage{latexsym,amssymb,amscd}
\usepackage{amsmath}
\usepackage[all]{xy}
\begin{document}
\title{ Layer lengths, torsion theories and the finitistic dimension}
\thanks{The authors thank the financial supports received from Proyecto PAPIIT-Universidad Nacional Aut\'onoma de M\'exico IN100810-3 MEXICO, Proyecto FCE-ANII 059 URUGUAY and NSERC Discovery Grants Program CANADA}
\author{ Fran\c cois Huard,\\ Marcelo Lanzilotta,\\ Octavio Mendoza}
\date{}

\begin{abstract} Let $\Lambda$ be a left-artinian ring. Generalizing the Loewy length, we propose the layer length associated with a torsion theory, which is a new measure for finitely generated $\Lambda$-modules. As an application, we obtain a theorem having as corollaries the main results of \cite{HLM2} and \cite{W}.
\end{abstract}

\maketitle
\section{Layer lengths}
Throughout the paper, we fix the following notation. $\Lambda$ will be a left-artinian ring and $\C:=\modu\,(\Lambda)$ the category of finitely generated left $\Lambda$-modules. We also denote by $\fun$ the category
of all additive functors from $\C$ to $\C$. Furthermore we let  $ \rad$ (resp.  $ \soc$) denote the Jacobson's radical (resp. socle) lying in $\fun$. Note that the functors $\rad$ and $\soc$ are both subfunctors of the identity $\Id.$ Recall that if $\alpha$ and $\beta$ belong to $\fun$ and $\alpha$ is a subfunctor of $\beta,$ we have the quotient functor $\beta/\alpha\in\fun$ which is defined as follows:

\smallskip

(a) $(\beta/\alpha)(M):=\beta(M)/\alpha(M)$ for $M\in\C$, and

\smallskip

(b) $(\beta/\alpha)(f)\,(x+\alpha\,(M)):=\beta\,(f)\,(x)+\alpha\,(N)$ for a morphism $f:M\rightarrow N$ in $\C.$

\smallskip

\noindent Furthermore, we set $\tops:=\Id/\rad\in\fun.$ Finally, we also recall that the functors $\rad$ and $\Id/\soc$ preserve monomorphisms and epimorphisms in $\C.$
\vspace{.2cm}

Given $\alpha\in\fun,$ we consider the {\bf $\alpha$-radical functor} $\Fa:=\rad\circ\alpha$ and the {\bf $\alpha$-socle quotient functor} $\Ga:=\alpha/(\soc\circ\alpha)$ where $\circ$ is the composition in $\fun.$ Furthermore, we consider the classes $\FFa=\{\ M\in \C\ :\ \alpha(M)=0\ \}$ and  $\Ta=\{\ M\in \C\ :\ \alpha(M)=M\ \}.$ Also we set $\mathrm{min}\,\emptyset:=\infty.$

\begin{defi}\label{layer-length} For $\alpha$ and $\beta$ in $\fun,$ we define:
 \begin{itemize}
  \item[(a)] the {\bf $\ab$-layer length} $\dl_{\alpha}^{\beta}:\C\rightarrow \Bbb{N}\cup\{\infty\}$
$$\dl_{\alpha}^{\beta}(M):=\mini\,\{ i\geq 0\,:\, \alpha\circ\beta^i(M)=0\ \};$$
  \item[(b)] the {\bf $\alpha$-radical layer length} $\dl^{\alpha}:=\dl_{\alpha}^{\Fa}$ and the {\bf $\alpha$-socle layer length} $\dl_{\alpha}:=\dl_{\alpha}^{\Ga}.$
 \end{itemize}
\end{defi}

Note that $\dl^{\alpha}(M)$ and $\dl_{\alpha}(M)$ are finite for all $M$ in $\C$.

\begin{ex} The Loewy length is obtained by taking $\alpha=\Id$ in \ref{layer-length} (b). This yields  the usual radical layer length $\dl^\Id$ and socle layer length $\dl_\Id.$ In this case, it is well known that $\dl^\Id=\dl_\Id.$
\end{ex}

We have the following natural question: Under what conditions  can we compare $\dl^\alpha$ with $\dl_\beta$ for some $\alpha,\beta\in\fun$? In order to answer this question, we propose ``Comparison Theorems". We first need to develop some necessary theory about layer lengths.

\begin{lem}\label{primero} Let $\alpha$ and $\beta$ be in $\fun.$ Then
 \begin{itemize}
  \item[(a)] $\FFa=\{\ M\in \C\ :\ \dl_{\alpha}^{\beta}\,(M)=0\};$
\vspace{.2cm}
  \item[(b)] $\dl_{\alpha}^{\beta}\,(M)=\dl_{\alpha}^{\beta}\,(\beta\,(M)) + 1$ for any $M\in\C\setminus\FFa$;
\vspace{.2cm}
  \item[(c)] if $M\in \Ta$ and $M\neq 0$, then $\dl^{\alpha}\,(M)=\dl^{\alpha}\,(\rad\,(M))+1$ and
$\dl_{\alpha}\,(M)=\dl_{\alpha}\,(M/\soc\,(M))+1$.
 \end{itemize}
\end{lem}
\begin{dem} (a) It is easy to see that: $M\in\FFa\Leftrightarrow\alpha\,(M)=0\Leftrightarrow\dl_{\alpha}^{\beta}\,(M)=0.$
\

 (b) Take $M\in \C$ such that $\alpha\,(M)\neq 0$. In particular $\dl_{\alpha}^{\beta}\,(M)\geq 1$.
Now for $i\geq 1$, $\beta^i\,(M)=\beta^{i-1}\,(\beta\,(M))$ and then $\dl_{\alpha}^{\beta}\,(M)=\dl_{\alpha}^{\beta}\,(\beta\,(M))+1$.
\

(c) Take $0\neq M\in \Ta$. It is clear that $M\in\C\setminus\FFa$ because $\FFa\cap \Ta=\{0\}$.  Since $\alpha\,(M)=M$, it follows that $\Fa\,(M)=\rad\,(M)$ and $\Ga\,(M)=M/\soc\,(M)$. Hence, by (b), we get the result.
\end{dem}

\begin{lem}\label{segundo} Let $\alpha$ and $\beta$ be in $\fun.$ Then
 \begin{itemize}
  \item[(a)] $\dl_{\alpha}^{\beta}\,(M\oplus N)=\maxi\,\{\dl_{\alpha}^{\beta}(M),\dl_{\alpha}^{\beta}(N)\}$ for any $M, N\in\C;$
\vspace{.2cm}
  \item[(b)] if $\alpha$ and $\beta$ preserve epimorphisms, then
   \begin{itemize}
\vspace{.2cm}
    \item[(b1)]  $\dl_{\alpha}^{\beta}\,(M)\leq \dl_{\alpha}^{\beta}\,(L)$ for any epimorphism $L\rightarrow M$ in $\C,$
\vspace{.2cm}
    \item[(b2)] $\dl_{\alpha}^{\beta}(M)\leq \dl_{\alpha}^{\beta}({}_{\Lambda}\Lambda)$ for each $M\in \C;$
\vspace{.2cm}
   \end{itemize}
  \item[(c)] If $\alpha$ and $\beta$ preserve monomorphisms, then $\dl_{\alpha}^{\beta}\,(L)\leq \dl_{\alpha}^{\beta}\,(M)$
for any monomorphism $L\rightarrow M$ in $\C.$
\end{itemize}
\end{lem}
\begin{dem} (a) Consider $a:=\dl_{\alpha}^{\beta}\,(M\oplus N)$. If $a=\infty$, then $\alpha\circ\beta^i\,(M\oplus N)\neq 0$ for all $i\geq 0$. Let $b:=\maxi\,\{\dl_{\alpha}^{\beta}\,(M), \dl_{\alpha}^{\beta}\,(N)\}$. If $b<\infty$ then $\alpha\circ\beta^b(M)=\alpha\circ\beta^b(N)=0$ hence $\alpha\circ\beta^b(M\oplus N)=0$ a contradiction. If $a<\infty$, then $0=\alpha\circ\beta^a\,(M\oplus N)$ and since $\alpha\circ\beta^a$ is an additive functor, it follows that $n:=\maxi\{\dl_{\alpha}^{\beta}\,(M), \dl_{\alpha}^{\beta}\,(N)\}\leq a.$ On the other hand $0=\alpha\circ\beta^n(M)\oplus\alpha\circ\beta^n(N)=\alpha\circ\beta^n(M\oplus N)$ implies that $a\leq n$.
\

(b) Suppose that $\alpha$ and $\beta$ preserve epimorphisms.  Let $n:=\dl_{\alpha}^{\beta}\,(L)$ where $f:L\rightarrow M$ is an epimorphism in $\C.$ If $n=\infty$ , then (b1) is immediate, otherwise, since $\alpha\circ\beta^n(f)$ is an epimorphism and $\alpha\circ\beta^n(L)=0,$ we get that $\dl_{\alpha}^{\beta}\,(M)\leq n;$ proving (b1). Finally, (b2) follows from (b1) and (a) since $M$ is finitely generated.
\

The proof of (c) is similar to that of (b1).
\end{dem}

We recall the following definitions that can be found, for example, in \cite{ASS, St}.

\begin{defi} For $\alpha\in\fun$, we say that $\alpha$ is a {\bf pre-radical} if $\alpha$ is a subfunctor of $\Id.$ In case $\alpha$ is a pre-radical, we set $\qa:=\Id/\alpha\in\fun$. Furthermore, if $\alpha$ is a pre-radical satisfying $\alpha\circ\qa=0$ then $\alpha$ is a {\bf radical}.
\end{defi}

The following two lemmas will be useful for the Comparison Theorems.

\begin{lem}\label{St} Let $\alpha\in\fun$ be a pre-radical. Then
\begin{itemize}
\item[(a)] $\alpha$ preserves monomorphisms;
\vspace{.2cm}
\item[(b)] $\FFa$ is closed under submodules and finite coproducts;
\vspace{.2cm}
\item[(c)] $\Ta$ is closed under quotients and finite coproducts;
\vspace{.2cm}
\item[(d)] if $\alpha$ preserves epimorphisms then $\alpha$ is a radical.
\end{itemize}
\end{lem}
\begin{dem} Let $f:M\rightarrow N$ be a monomorphism. Consider the inclusions $i_M: \alpha\,(M)\rightarrow M$ and
$i_N: \alpha\,(N)\rightarrow N$. Thus $i_N\,\alpha(f)=f\,i_M$, and since $f$ is a monomorphism, then so is $\alpha(f)$, proving (a). For (b) and (c) see \cite[Ch. VI, Prop. 1.2]{St}. For (d) see \cite[Ch. VI, Ex. 5]{St}.
\end{dem}

\begin{lem}\label{elocho} Let $\alpha\in\fun$ be a pre-radical. The following conditions are equivalent.
\begin{itemize}
\item[(a)] The functor $\alpha$ is left exact.
\vspace{.2cm}
\item[(b)] The class $\Ta$ is closed under submodules and $\alpha^2=\alpha.$
\vspace{.2cm}
\item[(c)] The functor $\qa$ preserves monomorphisms.
\end{itemize}
\end{lem}

\begin{dem} (a)$\Leftrightarrow$(b) See in \cite[Ch. VI, Prop. 1.7]{St}.
\

 (a)$\Leftrightarrow$(c) Consider an exact sequence
$0\rightarrow M \stackrel{f}{\rightarrow} N \stackrel{g}{\rightarrow} X\rightarrow 0$ in $\C.$ We have the following exact and commutative diagram
$$\xymatrix{ {} & 0 \ar[d] & 0 \ar[d] {} & {}\\
 0 \ar[r] & \alpha\,(M) \ar[d]_{\alpha\,(f)} \ar[r]^{i_M} & M \ar[d]_f \ar[r] & \qa\,(M) \ar[d]_{\qa\,(f)} \ar[r] & 0\\
 0 \ar[r] & \alpha\,(N) \ar[d]_{\pi} \ar[r]^{i_N} & N \ar[d]_g \ar[r] & \qa\,(N) \ar[r] & 0\\
 {} &  \overline{X} \ar[d] \ar[r]^\mu & X \ar[d] & {} & {}\\
 {} & 0 & 0 & {} & {}       }$$ where $\overline{X}:=\alpha(N)/\Ima(\alpha\,(f))$ and $\pi(x):= x+\Ima(\alpha\,(f)).$ By the Snake's Lemma, we have that $\Ker\,(\qa\,(f))\simeq\Ker\,(\mu),$ and hence the item (c) is equivalent to saying that $\mu$ is a monomorphism.
\

Since $gf=0$, it follows that $\alpha\,(g)\,\alpha(f)=0$ and hence there exists a morphism
$\beta:\overline{X}\rightarrow \alpha\,(X)$ such that $\beta\,\pi=\alpha(g);$ therefore $\Ker(\beta)=\Ker(\alpha(g))/\Ima(\alpha\,(f)).$  Furthermore, we assert that the following diagram is commutative
$$\xymatrix{ {} & \alpha\,(X)\ar[dr]^{i_X} & {} \\
\overline{X} \ar[ur]^\beta \ar[rr]^\mu & & X.}$$
Indeed $\mu\,\pi=g\,i_N=i_X\,\alpha\,(g)=i_X\,\beta\,\pi$ and since $\pi$ is an epimorphism, it follows that $\mu=i_X\,\beta.$ Therefore, the proof follows from the following equivalences : $\alpha$ is left exact $\Leftrightarrow$ $\Ker\,(\alpha\,(g))=\Ima\,(\alpha\,(f))$ $\Leftrightarrow$ $\beta$ is a monomorphism $\Leftrightarrow$ $\mu$ is a monomorphism $\Leftrightarrow$ $\Ker\,(\qa\,(f))=0.$
\end{dem}

\section{Two comparison theorems for layer lengths}

In this section, we state and prove two comparison theorems for layer lengths. To do that, we need some propositions and lemmas as follows.

\begin{pro}\label{quotient} Let $\alpha\in\fun$ be a pre-radical which preserves epimorphisms. If $0\neq M\in \C$ and $\soc\,(M)\in\Ta$, then $M\not\in\FFa$ and $\dl^{\alpha}\,(M/\soc\,(M))+1=\dl^{\alpha}\,(M)$.
\end{pro}
\begin{dem}
Let $0\neq M\in\C$ such that $\soc\,(M)\in\Ta$. Note that $M\not\in \FFa$, otherwise, by \ref{St} (b), $\soc\,(M)\in\FFa$ and so $\soc\,(M)\in\Ta\cap\FFa=\{0\};$ contradicting that $M\neq 0.$\\
Since $M\not\in\FFa$, we have $m:=\dl^{\alpha}(M) >0$ (see \ref{primero} (a)). Take $N:=\alpha\circ\Fa^{\!\!\!\!m-1}(M)$ which is a submodule of $M$ since $\alpha$ and $\rad$ preserve monomorphisms.\\
We assert that $N$ is semisimple. Indeed, $\alpha\,(\rad\,(N))=\alpha\circ\Fa^{\!\!\!\!m}(M)=0$ and hence $\rad\,(N)\in\FFa$. Then, by \ref{St} (b), it follows that $\soc\,(\rad\,(N))\in \FFa$. On the other hand, there is an epimorphism $\soc\,(M)=\soc\,(N)\oplus X\rightarrow\soc\,(\rad\,(N))$ for some $X\in\C.$ Thus, by \ref{St} (c), we get $\soc\,(\rad\,(N))\in \Ta$ since $\soc\,(M)\in\Ta.$ Hence $\soc\,(\rad\,(N))\in\Ta\cap\FFa=\{0\}$ showing
that $\rad\,(N)=0;$ therefore $N$ is semisimple. Furthermore, since the pre-radicals $\alpha$ and $\rad$ preserve epimorphisms, then  so does the preradical $\alpha\circ\Fa^{\!\!\!\!m-1}$. Therefore, we get the following commutative diagram
$$\xymatrix{ M \ar[r]^{p_1} & M/\soc\,(M) \\
   \alpha\circ\Fa^{\!\!\!\!m-1}\,(M)=N \ar[u]^{i_1} \ar[r]_{p_2}  &  \alpha\circ\Fa^{\!\!\!\!m-1}\,(M/\soc\,(M))  \ar[u]_{i_2}  }$$
where $p_1$ is the canonical epimorphism, $p_2=\alpha\circ\Fa^{\!\!\!\!m-1}\,(p_1)$ and $i_1$ and $i_2$ are the inclusions.
Now, $N$ semisimple implies that $i_1(N)\subseteq\soc\,(M)$ and hence $i_2\,p_2=p_1\,i_1=0$, proving that $p_2=0$
since $i_2$ is a monomorphism. Then we get that $\alpha\circ\Fa^{\!\!\!\!m-1}(M/\soc\,(M))=0$ and hence
$n:=\dl^{\alpha}\,(M/\soc\,(M))\leq m-1$. Suppose that $n<m-1$. Consider the following commutative diagram
$$\xymatrix{ M \ar[r]^{p_1} & M/\soc\,(M) \\
   \alpha\circ\Fa^{\!\!\!\!n}\,(M) \ar[u]^{j_1} \ar[r]_{p_3}  &  0  \ar[u]_{j_2}  }$$
where $p_3=\alpha\circ\Fa^{\!\!\!\!n}\,(p_1)$ and $j_1$ and $j_2$ are the inclusion maps. Thus $p_1\,j_1=0$ and then $\alpha\circ\Fa^{\!\!\!\!n}\,(M)\subseteq\soc\,(M)$. So $\alpha\circ\Fa^{\!\!\!\!n+1}\,(M)=\alpha\,(\rad\,(\alpha\circ\Fa^{\!\!\!\!n}(M)))=0$ and hence $\dl^{\alpha}\,(M)\leq n+1<m$, a contradiction, proving that $n=m-1$.
\end{dem}

\begin{lem}\label{soc} Let $\alpha,\beta\in\fun$ be such that $\beta$ is a radical and $\alpha$ is a pre-radical. If $\FFa\subseteq\Tb$ then $\soc\,(\Ima\,(\qb))\subseteq\Ta\cap \FFb$.
\end{lem}
\begin{dem} Since $\beta\circ\qb=0$ it follows that $\Ima\,(\qb)\subseteq \FFb,$ therefore,
$\soc\,(\Ima\,(\qb))\subset\FFb$. It remains to prove that $\soc\,(\Ima\,(\qb))\subseteq\Ta$. Let $M\in\Ima\,(\qb)$ and take any simple $S$ in $\soc\,(M)$. We assert that $\Hom_\Lambda(-,S)|_{\Tb}=0;$ otherwise, there is an $X\in \Tb$ and an epimorphism $X\rightarrow S$ and so $S\in \Tb$. Hence $S\in \Tb\cap \FFb$ since $\soc\,(M)\in\FFb$. Then $S=0$, a contradiction, proving that $\Hom_\Lambda(-,S)|_{\Tb}=0$. In particular, using that $\FFa\subseteq\Tb$, we get $\Hom_\Lambda(-,S)|_{\FFa}=0$ and then $S\not\in\FFa$. Consider now the following exact sequence
$$0\rightarrow \alpha\,(S)\rightarrow S\rightarrow\qa\,(S)\rightarrow 0.$$ We have that $\alpha\,(S)\neq 0$ (otherwise, $S\in\FFa)$. Hence $\alpha\,(S)=S$ giving us that $S\in \Ta$ for all simple in $\soc\,(M)$. Then $\soc\,(M)\in \Ta$.
\end{dem}

The following result will be our first comparison theorem.

\begin{teo}\label{laseis} Let $\alpha,\beta\in\fun$ be radicals such that $\alpha$ preserves epimorphisms. Then  $\dl_{\qb}\,(M)\leq \dl^{\alpha}\,(M)$ for each $M\in\C$ if and only if $\FFa\subseteq\Tb.$
\end{teo}
\begin{dem} $(\Rightarrow)$ Let $\dl_{\qb}\,(M)\leq \dl^{\alpha}\,(M)$ for each $M\in\C.$ If $M\in \FFa$ then $\dl^{\alpha}\,(M)=0$ and so $\dl_{\qb}\,(M)=0;$ proving that $M\in\FFq=\Tb.$ Hence $\FFa\subseteq\Tb.$
\

$(\Leftarrow)$ Let $\FFa\subseteq\Tb.$ To prove that $\dl_{\qb}\,(M)\leq \dl^{\alpha}\,(M)$ for each $M\in\C,$ we use induction on $\dl^{\alpha}\,(M).$ If $\dl^{\alpha}\,(M)=0$ then $M\in \FFa\subseteq \Tb=\FFq$. Therefore
 $\dl_{\qb}\,(M)=0$.\\
Suppose that $\dl^{\alpha}\,(M)\geq 1$. We may assume that $\qb\,(M)\neq0$ (otherwise the proof is done).
We assert the following:
 \begin{itemize}
  \item[(i)] $\dl^{\alpha}\,(\qb\,(M))\leq\dl^{\alpha}\,(M);$
\vspace{.2cm}
  \item[(ii)] $\dl^{\alpha}(\qb\,(M)/\soc\,(\qb\,(M)))=\dl^{\alpha}\,(\qb\,(M))-1;$
\vspace{.2cm}
  \item[(iii)] $\dl_{\qb}\,(\qb\,(M)/\soc\,(\qb\,(M)))\leq \dl^{\alpha}\,(\qb\,(M)/\soc\,(\qb\,(M))).$
 \end{itemize} Indeed, (i) follows from the canonical epimorphism $M\rightarrow \qb\,(M)$ and \ref{segundo}. On the other hand, we already have that $\qb\,(M)\neq0;$ moreover, \ref{soc} implies $\soc\,(\qb\,(M))\in\Ta$, then (ii) follows from \ref{quotient}. To prove (iii), we use (i) and (ii) as follows:
$$\dl^{\alpha}\,(\qb\,(M)/\soc\,(\qb\,(M)))=\dl^{\alpha}\,(\qb\,(M))-1\leq\dl^{\alpha}\,(M)-1<\dl^{\alpha}\,(M);$$
thus (iii) follows from our induction hypothesis.
\

Now, since $\qb^2=\qb/(\beta\circ\qb)$ and  $\beta$ is a radical, we have that $\qb^2=\qb$ and $\qb\,(M)\in\Tq=\FFb$. Therefore it follows from the definition of $\dl_{\qb}$ that
$\dl_{\qb}\,(M)=\dl_{\qb}\,(\qb\,(M))$ and so from \ref{primero} (c), we get\\ $\dl_{\qb}\,(\qb\,(M))=\dl_{\qb}(\qb\,(M)/\soc\,(\qb\,(M)))+1\leq \dl^{\alpha}\,(\qb\,(M)/\soc\,(\qb\,(M)))+1=\dl^{\alpha}\,(\qb\,(M))\leq \dl^{\alpha}\,(M).$
\end{dem}

\begin{pro}\label{+uno} Let $\alpha\in\fun$ be a pre-radical such that $\qa$ preserves monomorphisms. If $0\neq M\in\C$ and $\tops\,(M)\in\FFa,$ then $M\not\in\Ta$ and $\dl_{\qa}\,(M)=\dl_{\qa}(\rad\,(M))+1.$
\end{pro}
\begin{dem} Let $0\neq M\in\C$ satisfying that $\tops\,(M)\in\FFa$. Suppose that $M\in\Ta$, then $\tops\,(M)\in\Ta$ and hence $\tops\,(M)\in\Ta\cap\FFa=\{0\}$, a contradiction. Hence $M\not\in\Ta$. We assert that for any $i\geq 0$ there is a natural map $\pi_i:\Id\rightarrow\qa\circ\Gq^{\!\!\!\!\! i}$ such that for each $X\in\C$, the map
$\pi_{i,X}:X\rightarrow\qa\circ\Gq^{\!\!\!\!\! i}\;(X)$ is an epimorphism. Indeed, for each $Z\in\C$ consider the following canonical quotients maps $p_Z:Z\rightarrow \qa\,(Z)$ and $\sigma_Z:\qa\,(Z)\rightarrow\Gq\,(Z).$  For each $X\in\C,$ we set $\pi_{0,X}:=p_X$ and we inductively define $\pi_{i,X}:=p_{\Gq^{\!\!\!\!\!\!\! i}\;\;(X)}\,\sigma_{\Gq^{\!\!\!\!\!\!\! i-1}\,(X)}\,\pi_{i-1,X}$ for $i\geq 1.$ It can be seen that the maps defined above satisfy the required conditions.
\

Take $m:=\dl_{\qa}\,(M)$. Observe that $m>0$ since $M\not\in\Ta=\FFqa$. We assert that
$N:=\qa\circ\Gq^{\!\!\!\!\! m-1}(M)$ is semisimple. Indeed, $$0=\qa\circ\Gq^{\!\!\!\!\! m}\,(M)= \qa\circ\Gq\,(\Gq^{\!\!\!\!\! m-1}(M))=\qa\,(N/\soc\,(N))$$ and hence $N/\soc\,(N)\in \Ta$ proving that
$\tops\,(N/\soc\,(N))\in \Ta$. On the other hand, the canonical epimorphisms
$M\rightarrow N$ and $N\rightarrow N/\soc\,(N)$ induce an epimorphism
$\tops\,(M)\rightarrow\tops\,(N/\soc\,(N))$ which is splitting, giving us a monomorphism
$\tops\,(N/\soc\,(N))\rightarrow\tops\,(M)$. Hence $\tops\,(N/\soc\,(N))\in \FFa$ since
$\tops\,(M)\in \FFa,$ and therefore, $\tops\,(N/\soc\,(N))\in\Ta\cap\FFa=\{0\}$. Then $N=\soc\,(N)$ which means that $N$ is semisimple. In what follows, we use the fact that the additive functor $\qa\circ\Gq^{\!\!\!\!\! i}\;$ preserves monomorphisms for any $i\geq 0$. Consider the following commutative diagram:
$$\xymatrix{ \rad\,(M) \ar[d]_{\pi_{m-1,\rad\,(M)}} \ar[r]^{i_M} & M \ar[d]^{\pi_{m-1,M}}\\
   \qa\circ\Gq^{\!\!\!\!\! m-1}\,(\rad\,(M)) \ar[r]_j & \qa\circ\Gq^{\!\!\!\!\! m-1}\,(M)=N}$$ where $j:=\qa\circ\Gq^{\!\!\!\!\! m-1}\,(i_M)$ is a monomorphism since the inclusion $i_M$ is so.
Then $\pi_{m-1,M}\,i_M\,(\rad\,(M))\subseteq\rad\,(N)=0$ and hence $j\,\pi_{m-1,\rad\,(M)}=0.$ Using that $j$ is a monomorphism, we have $\pi_{m-1,\rad\,(M)}=0$ and thus $\qa\,\Gq^{\!\!\!\!\!m-1}(\mathrm{rad}(M))=0;$ proving that $\ell:=\dl_\qa\,(\rad\,(M))\leq m-1.$ We assert that $\ell=m-1.$\\
Suppose that $\ell<m-1$. In this case, we have the following commutative diagram:
$$\xymatrix{ \rad\,(M) \ar[d]_{\pi_{\ell,\rad\,(M)}} \ar[r]^{i_M} & M \ar[d]^{\pi_{\ell,M}}\\
   0 \ar[r] & \qa\circ\Gq^{\!\!\!\! \ell}\,(M).}$$ So, there is an epimorphism $\theta:\tops\,(M)\rightarrow \qa\circ\Ga^{\!\!\!\!\ell}\,(M)$ such that $\theta\, p=\pi_{\ell,M},$ where $p:M\rightarrow \tops\,(M)$ is the canonical epimorphism. Therefore $\qa\circ\Ga^{\!\!\!\! \ell}\,(M)$ is semisimple and hence $\dl_{\qa}\,(M)\leq \ell+1<m$ a contradiction; proving that $\ell=m-1.$
\end{dem}

\begin{lem}\label{top} Let $\alpha, \beta\in\fun$ be pre-radicals. If $\alpha^2=\alpha$ and $\Tb\subseteq\FFa$ then $\tops\,(\Ima\,(\alpha))\subseteq \FFb\cap \Ta$.
\end{lem}
\begin{dem}
Since $\alpha^2=\alpha$, we have $\Ima\,(\alpha)\subseteq \Ta$ and so, by \ref{St} (c),
$\tops\,(\Ima\,(\alpha))\subseteq \Ta$. Let $S$ be a simple $\Lambda$-module in $\tops\,(\alpha\,(M))$ for $M\in \C$. In particular $S\in\Ta$. Since $\Hom_\Lambda(S,-)|_{\FFa}=0$ (otherwise we would have that $S\in \FFa\cap\Ta=\{0\})$, it follows from $\Tb\subseteq\FFa$ that $\Hom_\Lambda(S,-)|_{\Tb}=0$ and hence $S\not\in\Tb$. Consider the exact sequence $0\rightarrow \beta\,(S)\rightarrow S\rightarrow \qb\,(S)\rightarrow 0$. If $\beta\,(S)\neq 0$ then $\beta\,(S)=S$ and this implies that $S\in \Tb$, a contradiction, proving that $\beta\,(S)=0.$ Hence $S\in \FFb$ for any simple $S$ in $\tops\,(\alpha(M)),$ and so, $\tops\,(\alpha\,(M))\in \FFb$ for any $M\in\C.$
\end{dem}

We are now ready to prove the second comparison theorem.

\begin{teo}\label{ladiez} Let $\alpha,\beta\in\fun$ be pre-radicals such that $\alpha^2=\alpha$ and $\qb$ preserves monomorphisms. Then $$\dl^{\alpha}\,(M)\leq \dl_{\qb}\,(M)\text{ for each }M\in\C\text{ if and only if }\Tb\subseteq\FFa.$$
\end{teo}
\begin{dem}
($\Rightarrow$) Let $M\in\Tb=\FFq.$ By \ref{primero} we have $\dl_{\qb}\,(M)=0$, and so, $\dl^{\alpha}\,(M)=0;$
giving us that $M\in\FFa$.
\

($\Leftarrow$) $\Tb\subseteq\FFa.$ The proof will be carried out by induction on $\dl_{\qb}\,(M)$ for $M\in\C$. If $\dl_{\qb}\,(M)=0$ then
$M\in\FFq=\Tb\subseteq\FFa$ and hence $\dl^{\alpha}\,(M)=0$ by \ref{primero}.\\
Suppose that $\dl_{\qb}\,(M)\geq 1$. We may assume that $\alpha(M)\,\neq 0$ (otherwise there is nothing to prove).
We assert the following:
 \begin{itemize}
  \item[(i)] $\dl_{\qb}\,(\alpha\,(M))\leq \dl_{\qb}\,(M);$
\vspace{.2cm}
  \item[(ii)] $\dl_{\qb}\,(\rad\,(\alpha\,(M))=\dl_{\qb}\,(\alpha\,(M))-1;$
\vspace{.2cm}
  \item[(iii)] $\dl^{\alpha}\,(\rad\,(\alpha\,(M)))\leq\dl_\qb\,(\rad\,(\alpha(M)));$
\vspace{.2cm}
  \item[(iv)] $\dl^{\alpha}\,(\rad\,(\alpha\,(M)))+1\leq \dl_{\qb}\,(M).$
 \end{itemize} Indeed, (i) follows from the monomorphism $\alpha\,(M)\rightarrow M$ and the fact that $\qb$ preserves monomorphisms (see in \ref{segundo} (c)). To prove (ii), we have  by \ref{top} that $\tops\,(\alpha\,(M))\in \FFb$. The result then follows from \ref{+uno}. On the other hand, we use (i), (ii) and induction to get (iii) as follows:
$\dl_{\qb}\,(\rad\,(\alpha\,(M)))\leq \dl_{\qb}\,(M)-1<\dl_{\qb}\,(M).$ Finally, the preceding inequalities and (iii) give us (iv).
\

To finish the proof, we use (iv) and $\alpha^2=\alpha$, as follows:
$$\dl^{\alpha}\,(M)=\dl^{\alpha}\,(\alpha(M))=\dl^{\alpha}\,(\rad\,(\alpha\,(M)))+1\leq \dl_{\qb}\,(M).$$
\end{dem}

\section{layer lengths induced by torsion radicals}

In order to state and prove our third comparison theorem, we recall some notions and basic results about torsion theories for $\C$ (see, for example, in \cite{ASS} and \cite{St}). In particular, the torsion radical attached to a torsion theory for $\C$ will be of crucial importance. Recall also that for a given class $\X$ in $\C$ and a $\Lambda$-module $M,$ the trace of $\X$ in $M$ is the $\Lambda$-submodule $\Trace_\X(M)$ of $M$ generated by the images $\mathrm{Im}\,(f)$ of the maps $f\in\Hom_\Lambda(X,M)$ with $X\in\X.$

\begin{defi} A {\bf torsion theory} for $\C$ is a pair $\tp$ of classes of modules in $\C$ satisfying the following conditions:
 \begin{itemize}
  \item[(a)] $\Hom_\Lambda\,(M,N)=0$ for any $M\in\T$ and $N\in\F;$
\vspace{.2cm}
  \item[(b)] for any $X\in\C,$ if $\Hom_\Lambda\,(X,-)|_\F=0$ then $X\in\T;$
\vspace{.2cm}
  \item[(c)] for any $X\in\C,$ if $\Hom_\Lambda\,(-,X)|_\T=0$ then $X\in\F.$
 \end{itemize}
\end{defi}

\begin{rk} Let $\tp$ be a torsion theory for $\C.$ We recall that $t:=\Trace_\T$ is the so called {\bf torsion radical} attached to $\tp.$ It can be seen that $t\,(M)$ is the largest submodule of $M$ lying in $\T.$
\end{rk}

The connection between torsion theories for $\C$ and subfunctors of $\Id$ is given by the following well known result.

\begin{pro}\label{corresp}  \cite{ASS, St} The map $\phi\,(\T,\F):=\Trace_\T$ induces a bijection $$\{\text{ torsion theories for } \C\}\stackrel{\phi}{\rightarrow}\{\text{ idempotent radicals in }\fun\}$$ with inverse $\phi^{-1}\,(\alpha)=(\T_\alpha,\F_\alpha).$
\end{pro}

Now, we are ready to state and prove the third comparison theorem.

\begin{teo}\label{elunico} Let $\tp$ and $\tpp$ be torsion theories for $\C;$ and consider the torsion radicals $t:=\Trace_\T$ and $t':=\Trace_{\T'}.$ Then
 $$\dl^t\,(M)=\dl_{\qtp}\,(M)\text{ for each }M\in\C\text{ if and only if }\mathcal{F}=\mathcal{T'}.$$
\end{teo}
\begin{dem}
($\Rightarrow$) Let $\dl^t\,(M)=\dl_{\qtp}\,(M)$  for each $M\in\C.$ Then, by \ref{primero} (a) and \ref{corresp},  we have the equivalences: $M\in\F=\F_t$ $\Leftrightarrow$  $\dl^t\,(M)=0=\dl_{\qtp}\,(M)$ $\Leftrightarrow$ $M\in \mathcal{F}_{\qtp}=\mathcal{T'}$.

($\Leftarrow$) Let $\mathcal{F}=\mathcal{T'}.$ To prove the result, it is enough to check the hyphotesis needed in \ref{laseis} and \ref{ladiez} for $\alpha=t$ and $\beta=t'.$ That is, we assert the following:
 \begin{itemize}
  \item[(i)] $t$ and $t'$ are idempotent radicals;
\vspace{.2cm}
  \item[(ii)] $t$ preserves epimorphisms;
\vspace{.2cm}
  \item[(iii)] $\qtp$ preserves monomorphisms.

 \end{itemize} Indeed, (i) follows from \ref{corresp}. To prove (ii), we have firstly that $\mathcal{F}=\mathcal{T'}$ implies $\mathcal{F}$ is closed under quotients, and since $t$ is radical, using \cite[Ch. VI, Ex. 5]{St}, we get the result. Finally, we get (iii) from \ref{elocho} since $\mathcal{T'}=\F_t$ is closed under submodules (see \ref{St} (b)) and $t'^2=t'$.
\end{dem}
\

An interpretation of \ref{elunico} in terms of ttf-classes can be given. So, we get a new view of ttf-classes by using layer lengths. In what follows, we recall the notion of ttf-class.

\begin{defi}\cite{J} A class $\X$ in $\C$ is called a {\bf ttf-class} if there exists classes $\T$ and $\F$ such that $(\T,\X)$ and $(\X,\F)$ are torsion theories for $\C$. In this case, the triple $(\T,\X,\F)$ is called a {\bf ttf-theory}.
\end{defi}

Using the notion of ttf-triple, we can get from \ref{elunico} the following result.
\begin{cor}\label{elcoro} Given a ttf-theory $(\T,\X,\F)$ for $\C,$ we have that $$\dl^t\,(M)=\dl_{q_x}\,(M)\text{ for each }M\in\C$$
where $t:=\Trace_\T$ and $x:=\Trace_{\X}$ are the corresponding torsion radicals.
\end{cor}
\begin{dem} This follows easily from \ref{elunico}.
\end{dem}

Given a class $\mathcal{A}$ in $\C,$ we denote by $\mathfrak{F}\,(\mathcal{A})$ the class of the $\mathcal{A}$-filtered modules in $\C.$ That is, $M \in\mathfrak{F}(\mathcal{A})$ if there is a finite
chain $0=M_0\subseteq M_1\subseteq \cdots\subseteq M_m=M$ of submodules of $M$ such that each quotient $M_i/M_{i-1}$ is
isomorphic to some object in $\mathcal{A}.$ We have the following easy characterization of ttf-classes for $\C.$

\begin{lem}\label{ttf-lema} Let $\X$ be a class in $\C.$ Then, the following conditions are equivalent:
 \begin{itemize}
  \item[(a)] $\X$ is a ttf-class for $\C;$
\vspace{.2cm}
  \item[(b)] $\X$ is closed under extensions, submodules and quotients;
\vspace{.2cm}
  \item[(c)] $\X=\mathfrak{F}(\mathcal{S})$ for some $\mathcal{S}\subseteq\{\ \Lambda$-simple modules $\}$.
 \end{itemize}
\end{lem}
\begin{dem} The equivalence of (b) and (c) is easy to see. On the other hand, The equivalence of (a) and (b) is well known (see for example in \cite{ASS} and \cite{St}).
\end{dem}

\begin{exs} $(1)$ Consider the trivial ttf-theory $(\C,0,\C)$.In this case, the torsion radicals are $t=\Trace_{\C}=\Id$ and $x=\Trace_0=0.$ So, $q_{x}=\Id$ and then $\dl^t$ is the radical layer length and $\dl_{q_x}$ is the socle layer length. In this case, \ref{elcoro} gives us the well known equality for the Loewy Length.
\

$(2)$ Consider the functors $K,S\in\fun$, which were introduced in \cite{HLM2}. These functors were defined throughout  the exact sequences: $0\rightarrow K\,(M)\rightarrow M\rightarrow Q\,(M)\rightarrow 0$ and $0\rightarrow S\,(M)\rightarrow M\rightarrow C\,(M)\rightarrow 0$ where $K\,(M)$ (resp. $C\,(M)$) is the maximal submodule (resp. quotient) of $M$ lying in $\mathfrak{F}(\mathcal{S}^{<\infty})$ with $\mathcal{S}^{<\infty}$ the class of the simple $\Lambda$-modules of
finite projective dimension. So, we have that $K$ and $S$ are idempotent radicals, and moreover,
$q_K=Q$ and $q_S=C$. Since $\F_S=\mathfrak{F}(\mathcal{S}^{<\infty})=\T_K$, we have that $(\T_S,\mathfrak{F}(\mathcal{S}^{<\infty}),\F_K)$ is a ttf-theory for $\C.$ Hence, by \ref{elcoro}, we get that $$\dl^{\infty}\,(M):=\dl^S\,(M)=\dl_{q_K}\,(M)=:\dl_{\infty}\,(M),$$ where $\dl^{\infty}$ is the so called infinite-layer length in \cite{HLM2}.
\end{exs}

In what follows we use the following notation. Let $\mathcal{S}\subseteq\{\ \Lambda$-simple modules $\},$  $\mathcal{S'}:=\{\ \Lambda$-simple modules $\}\setminus\mathcal{S}$ and $\X_\mathcal{S}:=\mathfrak{F}(\mathcal{S}).$ Following \ref{ttf-lema}, we say that the ttf-theory $(\T_\mathcal{S}, \X_\mathcal{S}, \F_\mathcal{S})$ is generated by the class $\mathcal{S}.$ We also consider the torsion radicals $t_\mathcal{S}:=\Trace_{\T_\mathcal{S}}$ and $x_\mathcal{S}:=\Trace_{\X_\mathcal{S}}.$

\begin{pro}\label{Prop1} Let  $(\T_\mathcal{S}, \X_\mathcal{S}, \F_\mathcal{S})$ be the ttf-theory generated by the class $\mathcal{S}.$ Then, the following statements hold true.
\begin{itemize}
  \item[(a)] $\T_\mathcal{S}=\{M\in\C\;:\;\tops\,(M)\in\add\,(\mathcal{S'})\}.$
\vspace{.2cm}
  \item[(b)] $\F_\mathcal{S}=\{M\in\C\;:\;\soc\,(M)\in\add\,(\mathcal{S'})\}.$
\vspace{.2cm}
  \item[(c)] The torsion class $\T_\mathcal{S}$ is closed under projective covers of objects in $\T_\mathcal{S}.$
\vspace{.2cm}
  \item[(d)] The set $I_\mathcal{S}:=t_\mathcal{S}\,({}_\Lambda\Lambda)$ is an ideal of $\Lambda$ and $t_\mathcal{S}(M)=I_\mathcal{S}\,M$ for any $M\in\C.$
\end{itemize}
\end{pro}
\begin{dem} (a) Let $M\in\C$ be such that $\tops\,(M)\in\add\,(\mathcal{S'}).$ We prove that $\Hom_\Lambda(M,-)|_{\X_\mathcal{S}}$ is zero. Suppose that there is some non zero morphism $f:M\to X$ with $X\in\X_\mathcal{S}.$ We may assume that $f$ is surjective, thus $\tops\,(X)$ is a direct summand of $\tops\,(M)$ contradicting that $\mathcal{S}\cap\mathcal{S'}=\emptyset;$ proving that $M\in\T_\mathcal{S}.$
\

 Let $M\in\C$ be such that $\Hom_\Lambda(M,-)|_{\X_\mathcal{S}}$ is zero. If there were some $S\in\mathcal{S}$ being a direct summand of $\tops\,(M),$ then by composing the canonical projections  $M\to\tops\,(M)\to S,$ we would obtain that $\Hom_\Lambda(M,S)\neq 0.$ Hence $\tops\,(M)\in\add\,(\mathcal{S'}).$
\

(b) It is quite similar to the proof of (a).
\

(c) It follows from (a) since $\tops\,(P_0(M))\simeq\tops\,(M),$ where $P_0(M)$ is the projective cover of $M.$
\

(d) It follows by \ref{ttf-lema} and \cite[Ch. VI, Ex. 5]{St}.
\end{dem}

\begin{pro} \label{Prop2} Let $\mathcal{S}_1$ and $\mathcal{S}_2$ be subsets of the set $\{\ \Lambda$-simple modules $\}.$ Then
\[\mathcal{S}_1\subseteq\mathcal{S}_2\quad\text{ if and only if }\quad\dl^{t_{\mathcal{S}_2}}\,(M)\leq\dl^{t_{\mathcal{S}_1}}\,(M)\quad\forall M\in\C.\]
\end{pro}
\begin{dem} It is clear that $\mathcal{S}_1\subseteq\mathcal{S}_2$ is equivalent to the inclusion $\X_{\mathcal{S}_1}\subseteq\X_{\mathcal{S}_2}.$ On the other hand, we have that $\X_{\mathcal{S}_1}=\F_{t_{\mathcal{S}_1}}$ and $\X_{\mathcal{S}_2}=\T_{x_{\mathcal{S}_2}}$ (see \ref{corresp}). Hence, by \ref{laseis}, it follows that $\mathcal{S}_1\subseteq\mathcal{S}_2$ iff $\dl_{q_{x_{\mathcal{S}_2}}}\,(M)\leq\dl^{t_{\mathcal{S}_1}}\,(M)\quad\forall M\in\C.$ Finally, from \ref{elcoro}, we know that $\dl_{q_{x_{\mathcal{S}_2}}}\,(M)=\dl^{t_{\mathcal{S}_2}}\,(M)\quad\forall M\in\C;$ proving the result.
\end{dem}

\section{Applications to the finitistic dimension}

\begin{defi} Let $\X$ and $\Y$ be classes of $\C$ and $\mathcal{S}\subseteq\{\ \Lambda$-simple modules $\}.$ We introduce the following classes of $\Lambda$-modules:
 \begin{itemize}
 \item[$\bullet$] $\X\oplus\Y:=\{X\oplus Y\;:\;X\in\X,\,Y\in\Y\},$
 \vspace{.2cm}
 \item[$\bullet$] $\C^{\mathcal{S}}_\ell:=\{M\in\C\;:\;\dl^{t_{\mathcal{S}}}\,(M)\leq\ell\},$
 \vspace{.2cm}
 \item[$\bullet$] $\mathbb{T}^{\mathcal{S}}_\ell:=\C^{\mathcal{S}}_\ell\oplus\Omega\,(\C^{\mathcal{S}}_\ell),$
 where $\Omega(M)$ denotes the first syzygy of $M\in\modu\,(\Lambda)$ and $\Omega(\mathcal X)=\{\Omega(M) : M\in\mathcal X\}$,
 \vspace{.2cm}
 \item[$\bullet$] the class $\mathcal{S}^{<\infty}$ of the simple $\Lambda$-modules of
finite projective dimension.
\end{itemize}

\end{defi}

\begin{rk}\label{rkalg} Let $\mathcal{S}\subseteq\{\ \Lambda$-simple modules $\},$ $J:=\rad\,(\Lambda)$ and $I_\mathcal{S}:=t_\mathcal{S}\,({}_\Lambda\Lambda).$ Consider the quotient ring $\Gamma:=\Lambda/J_\ell(\mathcal{S})$ where $J_\ell(\mathcal{S}):=I_\mathcal{S}(JI_\mathcal{S})^\ell\unlhd\Lambda.$ So, by \ref{Prop1} (d), it is not hard to see that \[\C^{\mathcal{S}}_\ell=\{M\in\C\;:\;J_\ell(\mathcal{S})M=0\}\simeq\mathrm{mod}\,(\Gamma).\]
\end{rk}

\begin{lem}\label{cuatrop} Let $\mathcal{S}\subseteq\{\ \Lambda\text{-simple modules }\}$ and $M\in\C.$ If $t_\mathcal{S}\,(M)\neq 0$ then \[\dl^{t_{\mathcal{S}}}\,(\Omega\,t_{\mathcal{S}}(M))\leq\dl^{t_{\mathcal{S}}}\,({}_\Lambda\Lambda)-1.\]
\end{lem}
\begin{dem} Assume that $t_\mathcal{S}\,(M)\neq 0.$ Consider the following exact sequence $0\to \Omega\,t_{\mathcal{S}}(M)\to P\to t_{\mathcal{S}}(M)\to 0$ where $P$ is the projective cover of $t_{\mathcal{S}}(M).$ Hence, by \ref{Prop1}, it follows that $0\neq P\in\T_{t_{\mathcal{S}}}.$ So, by using \ref{primero} and \ref{segundo}, we get \[\dl^{t_{\mathcal{S}}}\,(\Omega\,t_{\mathcal{S}}(M))\leq\dl^{t_{\mathcal{S}}}\,(\rad\,(P))=\dl^{t_{\mathcal{S}}}\,(P)-1\leq\dl^{t_{\mathcal{S}}}\,({}_\Lambda\Lambda)-1.\]
\end{dem}

In what follows, we use the function $\Psi:\modu\,(\Lambda)\to \mathbb N$ defined by Igusa and Todorov in \cite{IT}.  We refer to \cite{IT,HLM1} for the definition and main properties of this function.  Given $\mathcal X\subseteq \modu\,(\Lambda)$, we set $\Psi\mathrm{dim}\,(\mathcal X)=\sup\{\Psi(M) : M\in \mathcal X\}$.

\begin{teo}\label{bigteo} Let $\mathcal{S}\subseteq\mathcal{S}^{<\infty}.$ If $\dl^{t_{\mathcal{S}}}\,({}_\Lambda\Lambda)\leq 2\ell+1$ and $\Psi\mathrm{dim}\,(\mathbb{T}^{\mathcal{S}}_\ell)$ is finite, then \[\mathrm{fin.dim.}\,(\Lambda)\leq\mathrm{max}\,\{\pd\,(\mathcal{S}),\,2+\Psi\mathrm{dim}\,(\mathbb{T}^{\mathcal{S}}_\ell)\}<\infty.\]
\end{teo}
\begin{dem} Let $\beta:=\pd\,(\mathcal{S})$ and $M\in\C$ of finite projective dimension. Assume that $\dl^{t_{\mathcal{S}}}\,({}_\Lambda\Lambda)\leq 2\ell+1$ and $\Psi\mathrm{dim}\,(\mathbb{T}^{\mathcal{S}}_\ell)<\infty.$ Consider the canonical exact sequence $0\to t_\mathcal{S}\,(M)\to M\to q_{t_\mathcal{S}}\,(M)\to 0.$ Since $q_{t_\mathcal{S}}\,(M)\in\X_\mathcal{S}$ and $\pd\,(M)<\infty,$ it follows that $\pd\,(t_\mathcal{S}\,(M))<\infty$ and also that $$\pd\,(M)\leq\mathrm{max}\,\{\beta,\pd\,(t_\mathcal{S}\,(M))\}\leq\mathrm{max}\,\{\beta,1+\pd\,(\Omega\,t_\mathcal{S}\,(M))\}.$$ If $t_\mathcal{S}\,(M)=0$ then $M\in\X_\mathcal{S},$ and so $\pd\,(M)\leq\beta$ getting the result in this case. So, we can assume that $t_\mathcal{S}\,(M)\neq0,$ and hence by \ref{cuatrop} $\dl^{t_{\mathcal{S}}}\,(\Omega\,t_\mathcal{S}\,(M))\leq 2\ell.$ Therefore $N:=F_{t_\mathcal{S}}^\ell(\Omega\,t_\mathcal{S}\,(M))\in\C^{\mathcal{S}}_\ell.$ Recall, as can be seen in \ref{rkalg}, that $N=J_\ell(\
 mathcal{S})\Omega\,t_\mathcal{S}\,(M),$ and then we have the exact sequence $0\to N\to \Omega\,t_\mathcal{S}\,(M)\to N'\to 0$ where $N':=\Omega\,t_\mathcal{S}\,(M)/N\in\C^{\mathcal{S}}_\ell.$ Using \cite[Remark 5]{IT}, we get $$\pd\,(\Omega\,t_\mathcal{S}\,(M))\leq1+\Psi(N\oplus\Omega(N'))\leq 1+\Psi\mathrm{dim}\,(\mathbb{T}^{\mathcal{S}}_\ell);$$ proving the result.
\end{dem}

\begin{rk} By taking $\mathcal{S}=\emptyset$ in \ref{bigteo}, we obtain as a consequence the main result in \cite{W} (see Theorem 3) : If $J^{2\ell+1}=0$ and $\Gamma:=\Lambda/J^\ell$ is of finite representation type, then $\mathrm{fin.dim.}\,(\Lambda)$ is finite. Indeed, since $\mathcal{S}$ is empty, we have that $t_{\mathcal{S}}=\Id,$ and so by \ref{rkalg} $\C^{\mathcal{S}}_\ell\simeq\mathrm{mod}\,(\Gamma).$ Therefore the test class $\mathbb{T}^{\mathcal{S}}_\ell$ is of finite representation type, giving us that $\mathrm{fin.dim.}\,(\Lambda)\leq2+\Psi\mathrm{dim}\,(\mathbb{T}^{\mathcal{S}}_\ell)<\infty.$
\end{rk}

\begin{lem}\label{test1} Let $\mathcal{S}\subseteq\mathcal{S}^{<\infty},$ $\beta:=\pd\,(\mathcal{S})$ and $\Sigma':=\bigoplus_{X\in\mathcal{S'}}\,X.$ Then the following statements hold true.
\begin{itemize}
 \item[(a)] For any $M\in\C^{\mathcal{S}}_1,$ we have that $\Omega^{\beta+1}(M)\oplus P\simeq \Omega^{\beta+1}(M_{\mathcal{S'}})\oplus P'$ where $P$ and $P'$ are projective $\Lambda$-modules and $M_{\mathcal{S'}}\in\add\,(\mathcal{S'}).$
\vspace{.2cm}
 \item[(b)] $\Psi\mathrm{dim}\,(\mathbb{T}^{\mathcal{S}}_1)\leq 1+\beta+\Psi(\Omega^{\beta+1}(\Sigma')\oplus\Omega^{\beta+2}(\Sigma'))<\infty.$
\end{itemize}
\end{lem}
\begin{dem} (a) Let $M\in\C^{\mathcal{S}}_1.$ So $t_{\mathcal{S}}F_{t_{\mathcal{S}}}(M)=0$ and then $\rad\circ t_{\mathcal{S}}(M)=F_{t_{\mathcal{S}}}(M)\in\X_\mathcal{S}.$ Moreover $q_{t_{\mathcal{S}}}(M)\in\X_\mathcal{S},$ hence  $\pd\,(F_{t_{\mathcal{S}}}(M))\leq\beta$ and $\pd\,(q_{t_{\mathcal{S}}}(M))\leq\beta.$ Therefore, applying \cite[Lemma 3.6]{HLM1} to the exact sequences $0\to t_{\mathcal{S}}(M)\to M\to q_{t_{\mathcal{S}}}(M)\to 0$ and $0\to F_{t_{\mathcal{S}}}(M)\to t_{\mathcal{S}}(M)
\to \tops(t_{\mathcal{S}}(M))\to 0$ the result follows since $M_{\mathcal{S'}}:=\tops(t_{\mathcal{S}}(M))\in\add\,(\mathcal{S'})$ by \ref{Prop1} (a).
\

(b) Let $M=X\oplus\Omega(Y)$ with $X,Y\in\C^{\mathcal{S}}_1.$ Hence, by (a), we have that $\Omega^{\beta+1}(M)\simeq \Omega^{\beta+1}(X_{\mathcal{S'}})\oplus\Omega^{\beta+2}(Y_{\mathcal{S'}})\oplus P$ for some projective $\Lambda$-module $P$ and $X_{\mathcal{S'}}, Y_{\mathcal{S'}}\in\add\,(\mathcal{S'}).$ Applying \cite[Proposition 3.5]{HLM1} $1+\beta$ times to $M$ and then, by the preceding isomorphism, we get
$\Psi(M)\leq 1+\beta+\Psi(\Omega^{\beta+1}(M))\leq 1+\beta+\Psi(\Omega^{\beta+1}(\Sigma')\oplus\Omega^{\beta+2}(\Sigma'))<\infty.$
\end{dem}

\begin{cor}\label{Mhlm2} Let $\mathcal{S}\subseteq\mathcal{S}^{<\infty},$ $\beta:=\pd\,(\mathcal{S})$ and $\Sigma':=\bigoplus_{X\in\mathcal{S'}}\,X.$ If  $\dl^{t_{\mathcal{S}}}\,({}_\Lambda\Lambda)\leq 3$ then $\mathrm{fin.dim.}\,(\Lambda)\leq 3+\beta+\Psi(\Omega^{\beta+1}(\Sigma')\oplus\Omega^{\beta+2}(\Sigma'))<\infty.$
\end{cor}
\begin{dem} It follows from \ref{bigteo} and \ref{test1}.
\end{dem}
\begin{rk} About Corollary \ref{Mhlm2}.
\begin{itemize}
 \item[$1$] By taking $\mathcal{S}=\mathcal{S}^{<\infty},$ we get the main result in \cite{HLM2} (see Theorem 5.5). On the other hand, by \ref{Prop2}, we see that the ``strongest" version that can be obtained from \ref{Mhlm2} is precisely when $\mathcal{S}=\mathcal{S}^{<\infty}.$
 \item[$2$] If $\mathcal{S}$ is empty we get the well know result of ``radical cube equal to zero". That is, if $J^3=0$ then  $\mathrm{fin.dim.}\,(\Lambda)\leq 3+\Psi(\Omega(\Lambda/J)\oplus\Omega^2(\Lambda/J))<\infty$ where $J:=\rad\,(\Lambda).$
\end{itemize}
\end{rk}

\vskip3mm \noindent Fran\c cois Huard:\\ Department of mathematics,
Bishop's University,\\ Sherbrooke, Qu\'ebec, CANADA,  J1M1Z7.\\
{\tt fhuard@ubishops.ca}

\vskip3mm \noindent Marcelo Lanzilotta:\\ Instituto de Matem\'atica y Estad\'{i}stica Rafael Laguardia,\\
J. Herrera y Reissig 565, Facultad de Ingenier\'{i}a, Universidad de la Rep\'ublica. CP 11300, Montevideo, URUGUAY.\\
{\tt marclan@fing.edu.uy}

\vskip3mm \noindent Octavio Mendoza Hern\'andez:\\ Instituto de Matem\'aticas, Universidad Nacional Aut\'onoma de M\'exico.\\
Circuito Exterior, Ciudad Universitaria,
C.P. 04510, M\'exico, D.F. M\'EXICO.\\ {\tt omendoza@matem.unam.mx}

\end{document}